\documentclass[11pt]{article} 
\usepackage{amsfonts,euscript,theorem}
\newcommand{\no}{\noindent}
\newcommand{\bc}{\begin{center}}
\newcommand{\ec}{\end{center}}
\newcommand{\be}{\begin{equation}}
\newcommand{\ee}{\end{equation}}
\newcommand{\bea}{\begin{eqnarray*}} 
\newcommand{\eea}{\end{eqnarray*}}
\newcommand{\ran}{\ensuremath{\rangle}}
\newcommand{\lan}{\ensuremath{\langle}}
\newcommand{\rar}{\ensuremath{\rightarrow}}
\newcommand{\ii}{\ensuremath{\infty}}

\newcommand{\pd}{\ensuremath{\partial}} 
\newcommand{\pdp}{\ensuremath{\partial}^p} 
\newcommand{\pdx}{\ensuremath{\partial}_x} 
 
\newcommand{\h}[1]{\ensuremath{\hbox{ #1 }}}
 
\newcommand{\w}[1]{\ensuremath{\widetilde{#1}}}

\newcommand{\Om}{\ensuremath{\Omega}}

\newcommand{\ga}{\ensuremath{\Gamma}}
\newcommand{\ps}{\ensuremath{\psi}}
\newcommand{\ve}{\ensuremath{\varepsilon}}
\newcommand{\vp}{\ensuremath{\varphi}}
\newcommand{\pe}{\ensuremath{\vp _{\ve}}}

\newcommand{\Z}{\ensuremath{\mathbb{Z}}}
 \newcommand{\de}{\ensuremath{\delta}}
\newcommand{\N}{\ensuremath{\mathbb{N}}}
\newcommand{\NN}{\ensuremath{\mathbb{N}_0}}
\newcommand{\R}{\ensuremath{\mathbb{R}}}
\newcommand{\Rm}{\ensuremath{{\mathbb{R}}^m}}
\newcommand{\C}{\ensuremath{\mathbb{C}}}
\newcommand{\A}{\ensuremath{A_0(\R)}}
\newcommand{\Aq}{\ensuremath{A_q(\R)}}
\newcommand{\G}{\ensuremath{{\EuScript{G}}}} 
\newcommand{\GR}{\ensuremath{{\EuScript{G}}(\R)}} 
\newcommand{\Gm}{\ensuremath{{\EuScript{G}}(\Rm)}} 
\newcommand{\D}{\ensuremath{{\EuScript{D}}(\R)}} 
\newcommand{\DD}{\ensuremath{{\EuScript{D}}'(\R)}} 

\newtheorem{Th1}{Theorem}
\newtheorem{Th2}[Th1]{Theorem}
\newtheorem{Th3}[Th1]{Theorem}
\newtheorem{Cor1}{Corollary}
\newtheorem{Thp}[Th1]{Theorem}

\begin{document}
\sloppy 

\vspace*{4mm}
\bc 
\boldmath
{\large\textbf{ON BALANCED PRODUCTS OF THE DISTRIBUTIONS $x_{\pm}^{a}$  \\*[1mm] IN COLOMBEAU ALGEBRA  \GR}}\\*[5mm] {\large\textsc{B. P. Damyanov}} \\*[5mm]
{\textit{ Bulgarian Academy of Sciences, INRNE\,-\,Theory Group \\  72, Tzarigradsko Shosse, 1784 Sofia, BULGARIA }}
\unboldmath
\ec 

\vspace*{4mm}
\footnotesize
\parbox{5.75in}{\no Results on products of the distributions $x_{\pm}^{a}$ and $\de^{(p)}(x)$ with coinciding singularities are derived when the products are balanced so that they exist in the distribution space. These results follow the idea of a known result on distributional products published by Jan Mikusi\'nski in 1966. They are  obtained in Colombeau algebra \GR of generalized functions that contains the distribution space and the notion of 'association' in it allows obtaining results in terms of distributions.}

\vspace*{3mm}
\parbox{5.75in}{ {\textbf {Key Words : Distribution, Product, Colombeau Algebra}}}

\vspace*{6mm}
\normalsize 
\setlength{\baselineskip}{19pt}
\setlength{\parindent}{10mm}
\bc {\large \textsc{Introduction} }\ec

\vspace*{2mm}
\no In 1966, Jan Mikusi\'nski published  the following well-known result  \cite{mik}\,: 
\be x^{-1}\,.\,x^{-1} \,-\, \pi^{{}2}\,\de(x)\,.\, \de(x) \ = \ x^{-2}, \ \ x \in \R.  \label{mik} \ee
\no Although, neither of the products on the left-hand side here exists, their difference still has a correct meaning in the distribution space \DD. \ Formulas including balanced products of distributions with coinciding singularities can be found in the mathematical and physical literature. For balanced products of that kind, we used the name `products of Mikusi\'nski type' in a previous paper  \cite{ijpm}, where we derived a generalization of the basic Mikusi\'nski product (\ref{mik}) in Colombeau algebra of tempered generalized functions \,(see equation (\ref{fstep}) below). 

We recall that the Colombeau algebra \G, introduced first in \cite{colm}, has followed various constructions of differential algebras that include distributions, proposed by K\"{o}nig, Berg,  Antonevich and Radyno, Egorov, and other authors. Lately, the Colombeau algebra has become very popular since it has almost optimal properties, as long as the problem of multiplication of Schwartz distributions is concerned: \G \ is an associative differential algebra with the distributions linearly embedded in it, the multiplication being compatible with the products of  $C^\ii$-differentiable functions.  Moreover, the algebra \G \ has a relation of  'association' that is a faithful generalization of the equality of distributions. This notion enables obtaining results `on distributional level' that are in consistency with distribution theory.  

In this paper, we prove several results on balanced  products of the distributions $x_{\pm}^{a}, \ a \in\Om = \{a \in \R: \, a \ne -1, -2,\ldots\}$, as they are embedded in Colombeau algebra and the products admit associated distributions. The results extends those proved in \cite{zaa} regarding the distributional products $x_{+}^{a}\,.\,x_{-}^{b}$ \,for $a+b > -2$. Our intention here is to treat the singular cases when $a+b = p$ \,and $p$ takes negative integer values. So, we consider a balanced product of $x_{+}^{a}\,.\,x_{-}^{b}$ \,in the `boundary case' $a+b = -2$ \,and then proceed one step further dealing with the case $a+b=-3$. Another general result on multiplication of the distributions $x_{+}^{a}, x_{-}^{b}$ \,as a balanced product is proved when $a+b = -p$, for arbitrary $p\in\N$.  

\vspace*{2mm}
\bc {\large \textsc{1. Definitions and notation}} \ec 

\vspace*{2mm}

\no We first recall the basic definitions of Colombeau algebra \G\ on the real line \R. 

\vspace*{1mm}
\textit{Definition 1} --- If  \N \ stands for the natural numbers, denote $\NN = \N \cup \{0\}$, and $\de _{ij} = \{\,1$ \,if \,$i=j, \,= 0$ otherwise$\,\}, \,i,j\in\NN$. Then we put for arbitrary $q \in \NN$\,: 
\[\Aq = \{ \vp(x) \in \D: \int_{\R} x^{j}\,\vp (x)\,dx  = \de _{0j}, \ j = 0, 1,...,q \}.\]  
Denote also $\pe = \ve ^{-1}\vp (\ve ^{-1}x)$ \,for  $\vp \in \Aq, \ve > 0$, and \,$\check{\vp}(x) = \vp(-x)$. \ Now, let $\EuScript{E}\,[\R]$ be the algebra of functions $f(\vp , x): \A \times \R  \rar \C$ that are infinitely differentiable, by a fixed `parameter' \vp. Then, the generalized functions of Colombeau are elements of the quotient algebra 
\[ \G \equiv\GR = \EuScript{E}_{\mathrm{M}} [\R]\,/ \ \EuScript{I}\,[\R].\] 
Here $\EuScript{E}_{\mathrm{M}}[\R]$ \,is the subalgebra of `moderate' functions such that for each compact subset $K$ of \R \,and $p \in\N$ \,there is a $q\in \N$ \,such that, for each $\vp \in \Aq$,  
\[ \sup_{x \in K}\,|\pdp f(\pe, x)\,| = O(\ve^{-q}), \hbox{ as}  \ \ve \rar 0_+.\] 
The ideal $\EuScript{I}\,[\R]$  of $\EuScript{E}_{\mathrm{M}}[\R]$ consists of all functions such that for each compact $K \subset\R $ \ and any $p\in \N$ \,there is a $q\in\N$ \ such that, for every $r \geq q$ and $\vp \in A_r(\R)$, 
\[ \sup_{x \in K}\,| \pdp f(\pe, x)\,| = O(\ve^{r-q}), \hbox{ as} \ \ve \rar 0_+.\] 

\vspace*{1mm}
The algebra \G \ contains the distributions on \R, canonically embedded as a \C-vector subspace \,by the map  \vspace*{-1mm}
\[ i : \DD \rar \,\G : u \mapsto \w{u} = \{\,\w{u}(\vp, x) = (u*{\check{\vp}})(x): \,\vp\in\Aq \,\}. \]

\vspace*{1mm}
\textit{Definition 2} --- A generalized function $f\in\G$ \,is said to admit some $ u\in \DD$ \ as `associated distribution', denoted $f \approx u$, if for some representative $f(\vp_\ve, x)$  of $f$  and any $\ps (x)\in\D$ \,there is a $q\in\NN$ \,such that, for any $\vp (x)\in\Aq$,  \vspace*{-2mm}
\[\lim_{{}\ve \rar 0_+} \int_{\R} f(\vp_\ve , x) \ps (x)\,dx = \lan u, \ps\ran.\] 

\vspace*{1mm}
This definition is independent of the representatives and the association is a faithful generalization of the equality of distributions \cite{col84}. 

Then, by product of some distributions in the algebra \G, sometimes called `Colombeau product', is meant the product of their embeddings in \G, whenever the result admits an associated distribution.  

The following coherence result holds \cite[Proposition 10.3]{mob} : If the regularized model product (in the terminology of Kami\'nski) of two distributions exists, then their Colombeau product also exists and coincides with the former. On the other hand, in the general setting of Colombeau algebra \Gm \ \cite{col84} \ (when the parameter functions \vp \ are not defined as tensor products), as well as in the algebra \GR \ on the real line,  this assertion turns into an equivalence, according to a result by Jel\'{\i}nek \cite{jel}; cf. also a recent study by Boie  \cite{bo}.

\vspace*{1mm}
Now, let $\w{x^{\,-p}}$ and $\w{\de}^{(p)} (x) $ be the embeddings in \G \ of the distributions $x^{\,-p}$ and $\de^{(p)} (x), \,p \in \N $. The following balanced distributional product has been proved in  \cite{ijpm}, which generalizes the basic Mikusi\'{n}ski formula (\ref{mik}) for arbitrary $p,q \in \N$:
\begin{equation} \quad \w{ x^{\,-p}}\,.\,\w{ x^{\,-q}} - \pi^2\frac{(-1)^{p+q}}{(p-1)!\,(q-1)!}\,\w{\de}^{\,(p-1)}(x)\,.\,\w{\de}^{\,(q-1)}(x) \ \approx \ x^{\,-p-q} \ (x\in \R). \label{fstep}\end{equation} 

\vspace*{1mm}
We next introduce the following.

\vspace*{1mm}
\textit{Definition 3} --- Letting $\Om = \{a \in \R: \, a \ne -1, -2,\ldots\}$, we denote 
\[ \nu_+^{\,a}\equiv \nu_+^{\,a}(x) =  \{ \frac{x^a}{\ga(a+1)} \ \h{if} x > 0 , \quad  = 0 \ \h{if} x<0\}. \] 
\[ \nu_-^{\,a}\equiv \nu_-^{\,a}(x) =  \{  \frac{(-x)^a}{\ga(a+1)} \ \h{if} x < 0,  \quad = 0 \ \h{if} x >0 \}. \]
Since $x \mapsto \,\nu_{\pm}^{\,a} (x)$ is a locally-integrable function for $a > -1$,  we can define the distributions $\nu_{\pm}^{\,a}$, for arbitrary $a \in \Om$, by setting 
\[ \nu_{+}^{\,a} = \pdx^r \,\nu_{+}^{\,a +r }(x), \qquad  \nu_{-}^{\,a} = (-1)^r \,\pdx^r \,\nu_{-}^{\,a +r }(x). \] 
Here $r\in \NN$  \,is such that $a+r > -1$ \,and the derivatives are in distributional sense. Note that, due to the equations \ $\pdx \,\nu_{\pm}^{\,a} = \pm \ \nu_{\pm}^{\,a-1}$ \,(with no number coefficients), the calculations with these `normed' distributions are made easier.

\vspace*{1mm}
We shall also make use of the definition and some basic properties of the binomial coefficients. Recall that \cite[\S \,21.5.1]{KK} \,if $x, y \in\R$ and $n, k\in\N$ \,then, by definition, 
\[ \!\!\left(\begin{array}{c}\!\!x\!\!\\ \!\!n\!\!\end{array}\right) = \left\{\!\frac{x(x-1)\cdots(x-n+1)}{n!}  \hbox{ if }n>0, \ = 1 \hbox{ if } n=0, \ = 0 \hbox{ if }n<0 \!\right\}\]
and it holds the following addition formula
\be  \left(\begin{array}{c}\!\!x+y\!\!\\ \!\!n\!\!\end{array}\right) = \sum_{k=0}^n  \left(\begin{array}{c}\!\!x\!\!\\ \!\!k\!\!\end{array}\right) \left(\begin{array}{c}\!\!y\!\!\\ \!\!n-k\!\!\end{array}\right) \quad (n>0),  \label{adth}\ee
as well as the equation
\be  \left(\begin{array}{c}\!\!-x\!\!\\ \!\!n\!\!\end{array}\right) = (-1)^n \left(\begin{array}{c}\!\!x+n-1\!\!\\ \!\!n\!\!\end{array}\right) \quad (x>0). \label{bico}\ee

\vspace*{2mm}
\bc {\large \textsc{2. Results on Mikusi\'nski type products of the distributions $x_{\pm}^{\,a}$}} \ec 

\vspace*{2mm}
\begin{Th1}$\!\!.$ For any $a, b \in \Om$ \,such that $a + b  > -2$, the embeddings in \G \ of the distributions $\nu_{\pm}^{\,a}(x)$ \ and $\de(x)$ satisfy 
\be \w{\nu_{+}^{\,a}}\, .\,\w{\nu_{-}^{\,b}} \,- \,\w{\nu_{-}^{\,a+b+1}}\,.\,\w{\de}(x)\approx 0 \qquad (x\in\R).  \label{th1}
\ee  
\end{Th1}

Now, we proceed to the singular distribution product $x_{+}^{a}\,.\,x_{-}^{b}$ \,in the `boundary case' \,$a+b = -2$, or else \,$b=-a-2$. We prove the following.

\begin{Th2}$\!\!.$ For any $a \in \R\backslash \Z$, the embeddings in \GR \ of the distributions $\nu_{\pm}^{\,a}(x)$ and $\de(x), \ \ x\in\R$, \,satisfy 
\be \w{\nu_{+}^{\,a}}\, .\, \w{\nu_{-}^{\,-a-2}} \ - \ \w{\de}(x)\,.\,\w{\de}(x) \ \approx \ -\ \frac{a\,+\,1}{2} \ \de'(x).  \label{th2}
\ee 
\end{Th2}

\vspace*{1mm}
\textsc{Proof} : \,For arbitrary $\ps(x) \in \D$ and $\vp \in \A$, we first get on the change \,$-x/\ve = u$ \,and applying the Taylor theorem
\begin{eqnarray} \lan \w{\de}(\pe, x)\w{\de}(\pe, x), \ps(x)\ran  & \! = & \!\frac{1}{\ve^2} \int_{-\ve l}^{\,\ve l} \vp^2\left(-\frac{x}{\ve}\right)\ps(x)\,dx \nonumber\\ &\! =\! &\frac{\ps(0)}{\ve}\!\int_{-l}^{\,l} \vp^2(u)\,du - \ps'(0)\!\int_{-l}^{\,l}  u\,\vp^2(u)\,du + O(\ve). \label{dsq}\end{eqnarray}
Without lost of generality, it is assumed here that \,supp $\!\vp(x) \subseteq [-l, l]$ \,for some $l\in \R_+$; then $-l\leq -x/\ve \leq l$ implies \ $-l\ve \leq x \leq l\ve$.

Suppose now that $r \in \N $ \ is subject to the condition $r > \mathrm{max } \,\{a+1, -a-1\}$. Then, Definition 3, the derivative in Colombeau algebra, and the substitution $u = (y-x)/{\ve}$ \ give for the embeddings in \GR\,:
\be \w{\nu_+^{\,a}} (\pe,x)  = \frac{(-1)^r\,\ve^{-r}} {\ga(a+r+1)} \int_{-x/\ve}^{\,l} (x+\ve u)^{a+r}\vp^{(r)}(u)\,du; \label{+a}\ee 
\[\w{\nu_-^{\,-a-2}}(\pe,x) = \frac{ \ve^{-r} }{\ga(-a+r-1)} \int_{-l}^{-x/\ve} (-x-\ve v)^{-a +r-2 }\vp^{(r)}(v)\,dv.\]

\no For arbitrary $\ps(x) \in \D$, denote the functional \vspace*{-2mm}
\[F = \lan\,\w{\w{\nu_+^{\,a}}(\pe,x)\ \nu_-^{\,-a-2}}(\pe,x), \,\ps(x)\ran  \quad \hbox{and} \ G = \ga(a+r+1)\ga(-a+r-1) F.\] 
Then, we obtain on the substitution $w=-x/\ve$\,:
\[G  =  \frac{(-1)^r}{\ve}\int_{-l}^{\,l} \ps(-\ve w) \!\int_w^{\,l}\vp^{(r)}(u) \!\int_{-l}^w (u-w)^{a+r}(w-v)^{-a+r-2}\vp^{(r)}(v)\,dv\,du\,dw. \]
Employing now Taylor theorem and changing twice the order of integration -- which is permissible here -- we get 
\bea  G & = &  \frac{(-1)^r\,\ps(0)}{\ve}\!\int_{-l}^{\,l}\vp^{(r)}(u)\!\int_{-l}^u \vp^{(r)}(v)\!\int_v^u (u-w)^{a+r}(w-v)^{-a+r-2}\,dw\,dv\,du \\ & &- (-1)^{r}\ps'(0)\!\int_{-l}^{\,l} \vp^{(r)}(u)\!\int_{-l}^u \vp^{(r)}(v)\!\int_v^u w (u-w)^{a+r}(w-v)^{-a+r-2}dw\,dv\,du \\ & & + \ O(\ve). \eea  
The substitution $w \,\rar \,t = (w-v)/(u-v)$, together with the relations \,$w-v = (u-v)t, \ u-w = (u-v)(1-t)$, gives  for the main part $G'$ of this expansion
\bea  G' & = &  \frac{(-1)^r\ps(0)}{\ve}\!\int_{-l}^l \vp^{(r)}(u)\!\int_{-l}^u (u-v)^{2r-1}\vp^{(r)}(v)\,dv\,du\!\int_0^1(1-t)^{a+r}\,t^{-a+r-2}dt \\*[1mm] & &\!\! - (-1)^r\ps'(0)\!\int_{-l}^l \vp^{(r)}(u)\!\int_{-l}^u(u-v)^{2r} \vp^{(r)}(v)\,dv\,du\!\int_0^1(1-t)^{a+r}\,t^{-a+r-1}dt  \\*[1mm] & &\!\! - (-1)^r\ps'(0)\!\int_{-l}^l\vp^{(r)}(u)\!\int_{-l}^u v(u-v)^{2r-1}\vp^{(r)}(v)dvdu\!\int_0^1(1-t)^{a+r}t^{-a+r-2}dt. \eea

\no By the definition of first-order Euler integral, 
\be \int_0^1(1-t)^{\,a}\,t^{\,b}\,dt = \frac{\ga(a+1)\ga(b+1)}{\ga(a+b+2)} \quad (a, b > -1). \label{Eul} \ee
\no Taking into account this, we obtain
\bea  F & = &  \frac{(-1)^r\,\ps(0)}{\ve\,(2r-1)!}r!\int_{-l}^l \vp^{(r)}(u) \!\int_{-l}^u (u-v)^{2r-1}\vp^{(r)}(v)\,dv\,du \\*[1mm] & & - \frac{(-1)^r\,(-a+r-1)\,\ps'(0)}{(2r)!}\int_{-l}^l \vp^{(r)}(u)\!\int_{-l}^u(u-v)^{2r} \vp^{(r)}(v)\,dv\,du \\*[1mm] & & - \frac{(-1)^r\,\ps'(0)}{(2r-1)!} \int_{-l}^l\vp^{(r)}(u)\int_{-l}^u \,v \,(u-v)^{2r-1}\,\vp^{(r)}(v)\,dv\,du  + O(\ve). \eea

\no Recall that, according to the Leibnitz rule, the following formula for differentiation of integrals holds \cite[\S \,6.4.1]{KK}\,:
\be \pd_u \left(\int_{-l}^u T(u,v)\,dv\right) = \left.\int_{-l}^u \pd_u \left(T(u,v)\right)\,dv  + T(u,v)\right|_{\,v\,=\,u}. \label{dif}\ee

\no Now, an $r$-time integration by parts in the variable $u$ -- the integrated parts being zero each time -- and equation (\ref{dif}), with $T(u,v)|_{v=u}= 0$ \,in this case, yield further
\bea  F  &= &\frac{\ps(0)}{\ve\,(r-1)!}\int_{-l}^{\,l} \vp(u) \int_{-l}^u (u-v)^{r-1} \vp^{(r)} (v)\,dv\,du \\ & & - \frac{(-a+r-1)\,\ps'(0)}{(r)!}\int_{-l}^{\,l} \vp(u)\!\int_{-l}^u(u-v)^{r} \vp^{(r)}(v)\,dv\,du \\ & & - \frac{\ps'(0)}{(r-1)!} \int_{-l}^{\,l}\vp(u)\int_{-l}^u \,v \,(u-v)^{r-1}\,\vp^{(r)}(v)\,dv\,du  + O(\ve). \eea

\no Replacing $v$ by $u - (u-v)$ in the last term, combining the summands, and integrating by parts in the variable $v$, we get
\be  F =  \frac{\ps(0)}{\ve}\!\int_{-l}^{\,l} \vp^2(u)\,du  - \ps'(0)\!\int_{-l}^{\,l}  u\,\vp^2(u)\,du + \frac{(a+1)}{2}\,\ps'(0) + O(\ve). \label{end}\ee

\no Taking into account equations (\ref{dsq}) and (\ref{end}), as well as linearity of the functional, we finally obtain
\[ \lim_{{}\ve \rar 0_+} \lan\,[\,\w{\nu_+^{\,a}}(\pe,x) \,\w{\nu_-^{\,-a-2}}(\pe,x), \,\ps(x)\ran -  \w{\de}(\pe, x)\w{\de}(\pe, x)], \ps(x)\,\ran = \frac{a\,+\,1}{2} \,\lan \,\de', \ps\,\ran. \]
This proves equation (\ref{th2}), according to Definition~2.

\vspace*{2mm}
\textit{Remark} --- In the notation $x_{\pm}^a$, equation (\ref{th2}) takes the following more complex form
\[ - \,\w{x_{+}^{\,a}}\, .\, \w{x_{-}^{\,-a-2}} \ + \ \ga(a\!-\!1) \ga(1\!-\!a) \,\w{\de}(x).\,\w{\de}(x) \ \approx \  \frac{\ga(a) \ga(1\!-\!a)}{2}\ \de'(x) \]
or else,
\be \w{x_{+}^{\,a}}\, .\,\w{x_{-}^{\,-a-2}} \ + \ \frac{\pi}{(a\!-\!1)\,\sin\pi a} \,\w{\de}(x)\,.\,\w{\de}(x) \ \approx \ \frac{\pi}{2\,\sin\pi a} \ \de'(x).  \label{x} \ee 

\vspace*{2mm}
Next we proceed one step further, dealing with the product $\nu_+^a\,.\,\nu_-^b$ in the singular case $a+b = -3$, or otherwise \,$b=-a-3$. 

\begin{Th3}$\!\!.$ 
 For any $a \in \R\backslash \Z$, the embeddings in \GR \ of the distributions $\nu_{\pm}^{\,a}(x)$ and $\de(x), \ \ x\in\R$, \,satisfy 
\be \w{\nu_+^{\,a}}\, .\,\w{\nu_-^{\,-a-3}} \ + \ (2a + 3)\ \w{\de}(x)\,.\,\w{\de'}(x) \ \approx \ \frac{1}{2}  \left(\begin{array}{c}\!\!a+2\!\!\\ \!\!2\!\!\end{array}\right)\ \de''(x).  \label{th3}
\ee 
\end{Th3}

\vspace*{1mm}
\textsc{Proof} : \,For arbitrary $\ps(x) \in \D$ and $\vp \in \A $, we get on the change \,$-x/\ve = u$ \,and Taylor theorem
\begin{eqnarray} \lan\, \w{\de}(\pe, x) \w{\de'}(\pe, x), \ps(x)\,\ran \!\! & = & \!\! -\frac{1}{\ve^3} \int_{-\ve l}^{\ve l} \vp\left(-\frac{x}{\ve}\right)\vp'\left(-\frac{x}{\ve}\right)\ps(x)\,dx  \nonumber \\ & \!\!= & \!\!- \frac{\ps'(0)}{2\,\ve}\!\int_{-l}^{\,l} \vp^2(u)du + \frac{\ps''(0)}{2}\!\int_{-l}^{\,l} u\vp^2(u)du  + O(\ve). \label{dd'}\end{eqnarray}

If $\ps(x) \in \D$, denote the functional 
\[ F = \lan\,\w{\nu_+^{\,a}}(\pe,x)\,\w{\nu_-^{\,-a-3}}(\pe,x), \,\ps(x)\ran  \quad \hbox{and} \ G = \ga(a+r+1)\ga(-a+r-2) F.\] 
Suppose that $r \in \N $ \,is subject to the condition \vspace*{1mm}$r > \mathrm{max } \,\{a+2, -a-1, 1\}$. Then we obtain on the substitution $w=-x/\ve$\,:
\[ G  = \frac{(-1)^r}{\ve^2}\int_{-l}^{\,l} \ps(-\ve w) \!\int_w^{\,l}\vp^{(r)}(u) \!\int_{-l}^w (u-w)^{a+r}(w-v)^{-a+r-3}\vp^{(r)}(v)\,dv\,du\,dw. \] 
Employing Taylor theorem and changing twice the order of integration, we get 
\bea  G & = & (-1)^r\sum_{n=0}^2 \frac{\ps^{(n)}(0)}{n!\,\ve^{2-n}}\int_{-l}^{\,l} \vp^{(r)}(u)\!\int_{-l}^u \vp^{(r)}(v)  \int_v^u w^n (u-w)^{a+r}(w-v)^{-a+r-3}\,dw\,dv\,du + O(\ve)\\*[1mm] & = &  \ \ga(a+r+1)\ga(-a+r-2) \left[\,F_1 + F_2 + F_3\,\right] \ + \ O(\ve). \eea  
Further, the substitution $w \,\rar \,t = (w-v)/(u-v)$, equation (\ref{Eul}), and  $r$-time integration by parts in the variable $u$ by equation (\ref{dif}), when applied successively to the terms $F_j, \,j=1,2,3$, \,will all yield\,:
\bea  F_1 & = & \frac{(-1)^r\,\ps(0)}{\ve^2\ (2r-2)!}\!\int_{-l}^{\,l} \vp^{(r)}(u) \!\int_{-l}^u (u-v)^{2r-2}\vp^{(r)}(v)\,dv\,du \\ & = & \frac{\ps(0)}{\ve^2\,(r-2)!}\int_{-l}^{\,l} \vp(u) \int_{-l}^u (u-v)^{r-2} \vp^{(r)} (v)\,dv\,du\\ & = &  \frac{\ps(0)}{\ve^2}\int_{-l}^{\,l} \vp(u)\,\vp'(u)\,du =  \left.\frac{\ps(0)}{\ve^2}\,\vp^2(u)\right|_{-l}^{\ l} \ = \ 0. \\*[5mm]  F_2 & = & - \frac{(-a+r-2)\,\ps'(0)}{\ve\,(r-1)!}\int_{-l}^l \vp(u) \!\int_{-l}^u (u-v)^{r-1}\vp^{(r)}(v) \,dv\,du \\ & & - \frac{\ps'(0)}{\ve\,(r-2)!}\int_{-l}^l \vp(u) \!\int_{-l}^u v\,(u-v)^{r-2}\vp^{(r)}(v) \,dv\,du.  \\*[5mm]  F_3 & = & \frac{(-a+r-2)(-a+r-1)\,\ps''(0)}{2\,r!}\int_{-l}^l \vp(u) \!\int_{-l}^u  (u-v)^{r}\,\vp^{(r)}(v)\,dv\,du \\ & & + \frac{(-a+r-2)\,\ps''(0)}{(r-1)!}\int_{-l}^l \vp(u) \!\int_{-l}^u v\,(u-v)^{r-1}\,\vp^{(r)}(v)\,dv\,du \\ & & + \frac{\ps''(0)}{2\,(r-2)!}\int_{-l}^l \vp(u) \!\int_{-l}^u  v^2\,(u-v)^{r-2}\,\vp^{(r)}(v)\,dv\,du. \eea

\no Replacing $v = u - (u-v)$, combining the summands, and integrating by parts in the variable $v$, we get 
\bea F_2 & =& (2a+3)\,\frac{\ps'(0)}{2\,\ve}\int_{-l}^l \vp^2(u)\,du.\\*[2mm]  F_3 & = &  \frac{1}{2}  \left(\begin{array}{c}\!\!a+2\!\!\\ \!\!2\!\!\end{array}\right)\ps''(0) - (2a\,+\,3)\,\frac{\ps''(0)}{2}\int_{-l}^l u\,\vp^2(u)\,du. \eea

Taking now account of equation (\ref{dd'}), we obtain for $ F = \sum_{j=1}^{\ 3} F_j+O(\ve)$
\[ F\,(\ve)  = - (2a\,-\,3) \lan\,\w{\de}(\pe, x)\,.\,\w{\de'}(\pe, x), \ps(x)\,\ran  + \frac{1}{2}  \left(\begin{array}{c}\!\!a+2\!\!\\ \!\!2\!\!\end{array}\right)\,\lan\,\de''(\pe, x), \ps(x)\,\ran + O(\ve). \]
Then, in view of linearity of the functional, we finally obtain
\[ \lim_{{}\ve \rar 0_+} \lan\,[\w{\nu_+^{\,a}}(\pe,x)\w{\nu_-^{\,-a\!-\!3}}(\pe,x) +  (2a+3)\,\w{\de}(\pe, x)\w{\de'}(\pe, x)], \ps\ran  = \frac{1}{2}\left(\begin{array}{c}\!\!a+2\!\!\\ \!\!2\!\!\end{array}\right)\,\lan\,\de''(x), \ps\ran. \]
By Definition~2, this proves equation (\ref{th3}).

\vspace*{2mm}
Direct consequences from the results of Theorems~2--3 are given by this.

\begin{Cor1}$\!\!.$ For arbitrary $a \in \R\backslash \Z$, the following M-type products hold for the embeddings in \GR \ of the distributions $\nu_{\pm}^{\,a}(x), \ x\in\R$\,: 

\be  - \,\w{\nu_{+}^{\,a}}\, .\,\w{\nu_{-}^{\,- a -2}}\ + \ \w{\nu_{-}^{\,a}}\, .\,\w{\nu_{+}^{\,-a-2}} \ \approx \ (a+1)\,\de'(x). \label{n2} \ee 

\be  \w{\nu_+^{\,a}}\, .\,\w{\nu_-^{\,-a-3}} \ + \ \w{\nu_{-}^{\,a}}\, .\,\w{\nu_{+}^{\,- a-3}} \ \approx \  \left(\begin{array}{c}\!\!a+2\!\!\\ \!\!2\!\!\end{array}\right)\ \de''(x) \label{n3} \ee 
\end{Cor1}

\vspace*{1mm}
\textsc{Proof} : \,Replacing $x$ by $-x$ in equation (\ref{th2}) and subtracting the latter equation from the result, we get (\ref{n2}). \ Equation (\ref{n3}) is obtained if we replace $x$ by $-x$ in equation (\ref{th3}) and then add the result to the latter equation.  

\vspace*{2mm}
Further extension of the results of Theorems~2--3 on the products $\nu_{+}^{\,a}\,.\,\nu_{-}^{\,b}$ \,to the cases $a+b = - 4, \cdots, -p$ \ for arbitrary $p\in\N$ \ is also possible (though more difficult) to prove, but the results are not Mikusi\'nski type products\,: one gets for the balancing term not single product any more but a sum of such products. 

Nevertheless, equations (\ref{n2}) and (\ref{n3}) can be extended to the general case $a+b= -p$ \,for any natural $p$. More exactly, it holds the following.

\no \begin{Thp}$\!\!.$  For arbitrary $a\in\R\backslash\Z$ and $p\in\N$, the embeddings in \GR \ of the distributions $\nu_{\pm}^{\,a}(x), \ x\in\R$, \,satisfy 
\be (-1)^p\  \w{\nu_+^{\,a}}\, .\,\w{\nu_-^{\,-a-\!p-1}} \ + \ \w{\nu_-^{\,a}}\, .\,\w{\nu_+^{\,-a-p-1}} \ \approx \ \ \left(\begin{array}{c}\!\!a+p\!\!\\ \!\!p\!\!\end{array}\right)\de^{(p)}(x).  \label{thp}
\ee 
\end{Thp}

\textsc{Proof} : \,(i) Choosing an $r \in \N$ \,such that $r > \mathrm{max } \,\{a+p, -a-1, p-1\}$, we denote for arbitrary $\ps(x) \in \D$ and $\vp \in \A $ \vspace*{-1mm}
\[F_1 = \lan\,\w{\nu_+^{\,a}}(\pe,x)\,\w{\nu_-^{\,-a\!-p\,-1}}(\pe,x), \,\ps(x)\ran  \ \ \hbox{and} \ \ G = \ga(a+r+1)\ga(-a+r-p) \,F_1.\] 
Taking account of representation (\ref{+a}), we get on the the substitution $w=-x/\ve$\,:
\[G  = \frac{(-1)^r}{\ve^{p}}\int_{-l}^{\,l}\ps(-\ve w) \int_w^{\,l}\vp^{(r)}(u) \!\int_{-l}^w (u-w)^{a+r}(w-v)^{-a+r-p-1}\vp^{(r)}(v)\,dv\,du\,dw. \]
The Taylor theorem, with the short-hand notation $\ps_0^n = \ps^{(n)}(0)$, and change of the order of integration, according to Dirichlet formula, now give
\bea G  = \sum_{n=0}^{p}\frac{(-1)^{r+n}\,\ps_0^n}{\ve^{p-n}\ n!}\!\int_{-l}^{\,l}&\!\!\!\!\!\!\vp^{(r)}(u)&\!\!\!\!\!\!\int_{-l}^u \vp^{(r)}(v) \int_v^u w^n(u-w)^{a+r}(w-v)^{-a+r-p-1}dw\,dv\,du + O(\ve).\eea 
Further, the substitution $w \rar \,t = (w-v)/(u-v)$, together with of the equation \,$w^n$ $ = \sum_{k=0}^{n}\left(\begin{array}{c}\!\!n\!\!\\ \!\!k\!\!\end{array}\right)(u-v)^k\,t^k\,v^{n-k}$, yields 
\bea G &=& \sum_{n=0}^{p}\frac{(-1)^{r+n}\,\ps_0^n }{\ve^{p-n}\ n!}\int_{-l}^{\,l} \vp^{(r)}(u) \int_{-l}^u \vp^{(r)}(v)\\ & & \times\sum_{k=0}^{n}\left(\begin{array}{c}\!\!n\!\!\\ \!\!k\!\!\end{array}\right)(u-v)^{2r-p+k}v^{n-k}\,dv\,du \int_0^1(1-t)^{a+r}t^{-a+r-p+k-1}\,dt + O(\ve).\eea
Taking into account equation (\ref{Eul}) and putting  $q = r-p \in \NN$, by the above choice of $r$, we obtain for the main part $F_1'$ \,of the asymptotic expansion of $F_1$\,:
\[ F_1'  = \sum_{n=0}^{p}\sum_{k=0}^{n}\frac{(-1)^{r+n+k}\,\ps_0^n\left(\begin{array}{c}\!\!a-q\!\!\\ \!\!k\!\!\end{array}\right)}{\ve^{p-n}\,(r+q+k)!\,(n-k)!}\int_{-l}^{\,l} \vp^{(r)}(u) \int_{-l}^u \vp^{(r)}(v) (u-v)^{r+q+k}\,v^{n-k}\,dv\,du. \]
We have used here that
\[  \frac{\ga(-a+q+k)}{\ga(-a+q)\,k!} = (-1)^k \left(\begin{array}{c}\!\!a-q\!\!\\ \!\!k\!\!\end{array}\right). \]
Now, a multiple integration by parts in the variables $u$ and $v$ -- each integrated term being zero -- and equation (\ref{dif}), with $T(u,v)|_{v=u}= 0$ \,in this particular case, will all give
\[ F_1' = \sum_{n=0}^{p}\sum_{k=0}^{n}\frac{(-1)^{n+k}\,\ps_0^n\left(\begin{array}{c}\!\!a-q\!\!\\ \!\!k\!\!\end{array}\right)}{\ve^{p-n}\,(q+k)!\,(n-k)!}\int_{-l}^{\,l}\vp(u) \int_{-l}^u \vp^{(p+q)}(v) (u-v)^{q+k}\,v^{n-k}\,dv\,du. \]

\no Replacing $v = u - (u-v)$ and $v^{n-k} = \sum_{j=0}^{n-k}(-1)^j\left(\begin{array}{c}\!\!n-k\!\!\\ \!\!j\!\!\end{array}\right)u^{n-k-j}\,(u-v)^j$,  we get
\bea F_1' &= & \sum_{n=0}^{p}\sum_{k=0}^{n}\sum_{j=0}^{n-k}\frac{(-1)^{n+k+j}\ps_0^n\left(\begin{array}{c}\!\!a-q\!\!\\ \!\!k\!\!\end{array}\right)}{\ve^{p-n}\,(n\!-\!k\!-\!j)!\,(q+k)!\,j!} \int_{-l}^{\,l} u^{n-k-j}\vp(u)\int_{-l}^u\!(u-v)^{q+k+j}\vp^{(p-k-j)}(v)\,dv\,du.\eea

\no We next put $t = k+j$, integrate $q+t+1$ times by parts the integral in $v$,  and change the order of summation, obtaining
\[F_1' = \sum_{n=0}^{p}\sum_{t=0}^{n}\frac{(-1)^{n+t}\,\ps_0^n}{\ve^{p-n}\,(n\!-\!t)!}\sum_{k=0}^{t} \left(\begin{array}{c}\!\!a-q\!\!\\ \!\!k\!\!\end{array}\right) \left(\begin{array}{c}\!\!q+t\!\!\\ \!\!t-k\!\!\end{array}\right) \int_{-l}^{\,l} \!u^{n-t}\vp(u)\vp^{(p-t-1)}(u)\,du. \]

\no Below, we shall use the short-hand notation 
\[  I_{n,\,t} := \int_{-l}^{\,l} u^{\{n-t\}}\vp(u)\vp^{(p-t-1)}(u)\,du, \ \hbox{ where }  u^{\{k\}} = u^k / k! \ \hbox{ and }\ k\in\NN,\]
noting that $\pd_u u^{\{k\}} = u^{\{k-1\}}$. Then, equations  (\ref{adth}) and (\ref{bico}) for binomial coefficients give
\be \qquad F_1 = \sum_{n=0}^{p}\frac{(-1)^{n}\,\ps_0^n}{\ve^{p-n}}\ \sum_{t=0}^{n}\left(\begin{array}{c}\!\!-a-1\!\!\\ \!\!t\!\!\end{array}\right) \int_{-l}^{\,l} u^{\{n-t\}}\vp(u)\vp^{(p-t-1)}(u)\,du  + O(\ve). \label{f1} \ee

\vspace*{2mm}
(ii) On the other hand, for arbitrary $\ps(x) \in \D$, denoting \vspace*{-1mm}
\[F_2  = \lan\,\w{\nu_-^{\,a}}(\pe,x)\,\w{\nu_+^{\,-a\!-p\,-1}}(\pe,x), \,\ps(x)\ran  \]
and proceeding similarly to the above calculations, we obtain for the main part of the latter functional 
\[ F_2' = \sum_{n=0}^{p}\sum_{k=0}^{n}\frac{(-1)^{n+k+1}\,\ps_0^n\left(\begin{array}{c}\!\!\!\!-a-r-1\!\!\!\!\\ \!\!k\!\!\end{array}\right)}{\ve^{p-n}(q+k)!\,(n-k)!}\int_{-l}^{\,l}\!\vp(u)\!\int_u^{\,l}\!\vp^{(p+q)}(v) (u-v)^{q+k}v^{n-k}dv\,du. \]
\no The replacements $v = u - (u-v)$ and $t = k+j$, integration $q+t+1$ times by parts in the variable $v$,  and change of the order of summation now give
\[ F_2' \! = \!\sum_{n=0}^{p}\sum_{t=0}^{n}\frac{(-1)^{n+t}\,\ps_0^n}{\ve^{p-n}}\sum_{k=0}^{t}\!\left(\begin{array}{c}\!\!\!\!-a\!-\!p\!-\!q\!-\!1\!\!\!\\ k\end{array}\right)\!\left(\begin{array}{c}\!\!q+t\!\!\\ \!\!t-k\!\!\end{array}\right)\!\!\!\int_{-l}^{\,l}\!u^{\{n-t\}}\vp(u)\vp^{(p-t-1)}(u)du. \]

\no Applying again equations  (\ref{adth}) and (\ref{bico}) to the binomial coefficients, we finally get
\be   \qquad F_2 =\sum_{n=0}^{p}\frac{(-1)^{n}\,\ps_0^n}{\ve^{p-n}}\ \sum_{t=0}^{n}\left(\begin{array}{c}\!\!a+p\!\!\\ \!\!t\!\!\end{array}\right) \int_{-l}^l u^{\{n-t\}}\vp(u)\vp^{(p-t-1)}(u)\,du  + O(\ve). \label{f2} \ee

\vspace*{2mm}
(iii) Next, we represent each of the functionals $F_i, i=1,2$ \,in the form
\be  \qquad  F_i = \sum_{n=0}^{p} \frac{(-1)^{n}\,\ps_0^n}{\ve^{p-n}}\ S_i (n)  + R_i + O(\ve), \ \ \hbox{where} \ \ R_i := \left.F_i \right|_{\,n = t = p}, \ i=1,2. \label{deco}\ee

\no Then, from equations (\ref{f1}) and (\ref{f2}) \,and taking into account that, by Definition~1, the parameter functions satisfy
\[  \int_{-1}^l \vp(u)\vp^{(-1)}(u)\,du  \equiv  \int_{-1}^l \vp(u)\int_{-l}^u\vp(u)\,dv\,du = \frac{1}{2}\left.\left(\int_{-l}^u \vp(v)\,dv\right)^2\right|_{-l}^{\ l} = \frac{1}{2}, \]
we obtain 
\[ R_1 =  \frac{(-1)^p\,\ps_0^p}{2} \left(\begin{array}{c}\!\!-a-1\!\!\\ \!\!p\!\!\end{array}\right) = \frac{1}{2} \left(\begin{array}{c}\!\!-a-1\!\!\\ \!\!p\!\!\end{array}\right) \lan\de^{(p)}(x), \,\ps(x)\ran. \]
 \[ R_2 = \frac{(-1)^p\,\ps_0^p}{2} \left(\begin{array}{c}\!a+p\!\!\\ \!\!p\!\!\end{array}\right)  = \frac{1}{2}\ \left(\begin{array}{c}\!\!a+p\!\!\\ \!\!p\!\!\end{array}\right)\lan \de^{(p)}(x), \,\ps(x)\ran. \]

\no Applying now equation (\ref{bico}), we get
\be  \qquad (-1)^p \ R_1 + R_2 = \left(\begin{array}{c}\!\!a+p\!\!\\ \!\!p\!\!\end{array}\right)\lan \de^{(p)}(x), \,\ps(x)\ran. \label{r1r2}\ee 

\vspace*{2mm}
(iv) It remains to prove that, for each $n = 0, 1, \ldots, p \ (t\leq p-1$ \,if $n=p)$, the following equation holds
\be (-1)^p \ S_1(n) \ + \ S_2(n) = 0. \label{s1+2}\ee 
To check this, we will distinguish between the cases when $p$ takes even and odd values, denoting them by $p=2h$ and $p=2h+1$ \,for some $h\in\NN$. The case $n=0$ \,is trivial. Passing to $n=1$, we first calculate the integrals involved. If  $\ p=2h$, we obtain on $h-1$ times integration by parts\,:
\[  I_{1,\,0} =  (-1)^{h-1}\,\left(\frac{1}{2} - h\right) \int_{-l}^{\,l} \left(\vp^{(h-1)}(u)\right)^2\,du, \ I_{1,\,1} =  (-1)^{h-1}\,\int_{-l}^{\,l} \left(\vp^{(h-1)}(u)\right)^2\,du. \]
We therefore have
\[ (-1)^{2h}S_1(1) + S_2(1)  = (-1)^{h-1} (1-2h -a-1+a+2h) \int_{-l}^{\,l} \left(\vp^{(h-1)}(u)\right)^2\,du = 0. \]
If \ $p = 2h+1$, then an $h$-time integration by parts gives
\[I_{1,\,0} = (-1)^{h} \int_{-l}^{\,l} u\,\left(\vp^{(h)}(u)\right)^2 \,du, \qquad I_{1,\,1} \ = \ 0. \]
Hence
\[ S_1(1)  = (-1)^{h}\ \int_{-l}^{\,l}\left(\vp^{(h)}(u)\right)^2\,du + 0 \ = \ S_2(1), \]
i.e. \ $(-1)^{2h+1}\,S_1(1) + S_2(1) = 0$, and equation (\ref{s1+2}) follows in the case $n=1$. 

\vspace*{1mm}
(v) This line of calculations is applicable for each $n$ \,and we shall perform it in the most complex case $n = p$. We compute first the integrals included, omitting some of the intermediate calculations. We make use of the next equation, which is obtained on $h-1$ times integration by parts and combining the terms
\be I_{p,\,t} = \sum_{m=0}^{h-1}(-1)^{h-1+m}\left(\begin{array}{c}\!\!h-1\!\!\\ \!\!m\!\!\end{array}\right)\int_{-l}^{\,l} u^{\{p-t-m\}}\,\vp^{(h-1-m)}(u)\vp^{(p-t-h)}(u)\,du. \label{ipt}\ee

Consider now the case $p=2h$. Taking account of (\ref{ipt}), we integrate further so as to obtain terms of the kind \ $\int\!u^{\{k\}}\left(\vp^{(n)}(u)\right)^2du$, for some $k, n\in\NN$; \,they will serve as `basic elements'. Proceeding this way, we get
\be I_{2h,\,t}  = \sum_{m=0}^{h-1}C_{h,\,m} \left[ \left(\begin{array}{c}\!\!m-h\!\!\\ \!\!2m+1-t\!\!\end{array}\right) +  \left(\begin{array}{c}\!\!m-h+1\!\!\\ \!\!2m+1-t\!\!\end{array}\right)\right],\label{i2h}\ee
where we have denoted
\[C_{h,\,m} = \sum_{m=0}^{h-1}\frac{(-1)^{h+m}}{2}\int_{-l}^{\,l} u^{\{2h-1-2m\}}\left(\vp^{(h-1-m)}\right)^2. \]
Inserting now equation (\ref{i2h}) in $S_1(p)$ for $p=2h$ and changing the order of summation, we get
\[ S_1(2h) = \sum_{m=0}^{h-1}C_{h,\,m} \sum_{t=0}^{2m+1}\left[ \left(\begin{array}{c}\!\!m-h\!\!\\ \!\!2m+1-t\!\!\end{array}\right) +  \left(\begin{array}{c}\!\!m-h+1\!\!\\ \!\!2m+1-t\!\!\end{array}\right)\right] \left(\begin{array}{c}\!-a-1\!\!\\ \!\!t\!\!\end{array}\right). \]
Then, the addition formula (\ref{adth}) and equation (\ref{bico}) for binomial coefficients yield
\bea S_1(2h) & = & \sum_{m=0}^{h-1}C_{h,\,m} \left[ \left(\begin{array}{c}\!\!m-a-h-1\!\!\\ \!\!2m+1\!\!\end{array}\right) +  \left(\begin{array}{c}\!\!m-a-h\!\!\\ \!\!2m+1\!\!\end{array}\right)\right] \\ & = & \sum_{m=0}^{h-1}  (-1)^{2m+1} \,C_{h,\,m} \left[ \left(\begin{array}{c}\!\!a+h+m+1\!\!\\ \!\!2m+1\!\!\end{array}\right) +  \left(\begin{array}{c}\!\!a+h+m\!\!\\ \!\!2m+1\!\!\end{array}\right)\right]. \eea
On the other hand, taking into account equations (\ref{i2h}) and (\ref{adth}), we obtain
\bea S_2(2h) &=& \sum_{m=0}^{h-1}C_{h,\,m} \sum_{t=0}^{2m+1}\left[ \left(\begin{array}{c}\!\!m-h\!\!\\ \!\!2m+1-t\!\!\end{array}\right) +  \left(\begin{array}{c}\!\!m-h+1\!\!\\ \!\!2m+1-t\!\!\end{array}\right)\right] \left(\begin{array}{c}\!a+2h\!\!\\ \!\!t\!\!\end{array}\right) \\ & = &  \sum_{m=0}^{h-1}C_{h,\,m} \left[ \left(\begin{array}{c}\!\!a+h+m\!\!\\ \!\!2m+1\!\!\end{array}\right) +  \left(\begin{array}{c}\!\!a+h+m+1\!\!\\ \!\!2m+1\!\!\end{array}\right)\right] \ = \ - \,S_1(2h).\eea
Thus 
\[ (-1)^{2h}S_1(2h) \  + S_2(2h) \ = \ 0. \]

Let further $p=2h+1$. Proceeding in a similar way, we obtain
\be I_{2h+1,\,t} = \sum_{m=0}^{h-1}C_{h,\,m} \left[ \left(\begin{array}{c}\!\!m-h\!\!\\ \!\!2m-t\!\!\end{array}\right) +  \left(\begin{array}{c}\!\!m-h+1\!\!\\ \!\!2m-t\!\!\end{array}\right)\right] .\label{bk}\ee 
Note that the `basic elements' included in $C_{h,\,m}$ \,coincide exactly with those in the case of even $p$. Then, we compute
\bea S_1(2h+1)\!\! & =& \!\!\sum_{m=0}^{h-1}C_{h,\,m}\!\sum_{t=0}^{2m}\left[ \left(\begin{array}{c}\!\!m-h\!\!\\ \!\!2m-t\!\!\end{array}\right)\!+\! \left(\begin{array}{c}\!\!m-h+1\!\!\\ \!\!2m-t\!\!\end{array}\right)\right]\!\left(\begin{array}{c}\!-a-1\!\!\\ \!\!t\!\!\end{array}\right) \\ & \!\!=\!\! &  \sum_{m=0}^{h-1} (-1)^{2m} C_{h,\,m} \left[ \left(\begin{array}{c}\!\!\!a+h+m+1\!\!\!\\ \!\!2m\!\!\end{array}\right)\!+\! \left(\begin{array}{c}\!\!a+h+m\!\!\\ \!\!2m+1\!\!\end{array}\right)\right] \!=\! S_2(2h+1).\eea
Therefore
\[ (-1)^{2h+1} \ S_1(2h+1) + S_2(2h+1) \ = \ 0. \]
This completes the proof of equation (\ref{s1+2}). Now, in view of equations (\ref{f1}) -- (\ref{s1+2}) and  linearity of the functional, we can finally write  
\bea \lim_{{}\ve \rar 0_+} \lan[(-1)^p\w{\nu_+^{a}}(\pe,x)\w{\nu_-^{-a\!-\!p\,-1}}(\pe,x) & + & \w{\nu_-^{a}}(\pe,x)\w{\nu_+^{-a\!-\!p\,-1}}(\pe,x)], \ps(x) \,\ran  \\ &  = & \left(\begin{array}{c}\!\!a+p\!\!\\ \!\!p\!\!\end{array}\right) \lan\de^{(p)}, \ps(x) \,\ran. \eea
According to Definition~2, this proves the result of Theorem~4. 

\vspace*{4mm}
\setlength{\baselineskip}{18pt}

\vspace*{2mm} 
\setlength{\baselineskip}{16pt}
\no {\footnotesize{Bulgarian Acad. Sci., INRNE - Theory Group \\
72 Tzarigradsko shosse, 1784 Sofia, Bulgaria \\ 
E-mail: damyanov@netel.bg}}

\end{document}